\documentclass[a4paper,english,reqno]{amsart}
\usepackage{babel}

\usepackage{units}

\usepackage[latin1]{inputenc}
\usepackage{algorithmic}
\usepackage{amsmath, amssymb}
\usepackage{amsthm}
\usepackage{xypic}
\xyoption{all}
\usepackage{url}
\def\@mpty{}
\providecommand{\N}{\mathbb{N}}
\providecommand{\Nplus}{\mathbb{N_+}}

\providecommand{\cardinality}[1]{|{#1}|}
\providecommand{\characteristic}[1]{\mathord{\mathrm{char}({#1})}}
\providecommand{\unitGroup}[1]{{#1}^*}

\providecommand{\divides}{\mathbin{\big\vert}}
\providecommand{\ndivides}{\mathbin{\not{\big\vert}}}

\providecommand{\gcd}[2]{}
\renewcommand{\gcd}[2]{({#1},{#2})}

\providecommand{\leadingCoeff}[2][]{\mathord{\mathrm{lc}_{#1}}({#2})}
\providecommand{\leadingCoeffLarge}[2][]{\mathord{\mathrm{lc}_{#1}}\left({#2}\right)}

\providecommand{\Automorphismgroup}[2]{\mathrm {Aut}({#1}/{#2})}
\providecommand{\AutomorphismGroup}[2]{\Automorphismgroup{#1}{#2}}

\providecommand{\PGL}[2]{\mathrm {PGL}_{#1}({#2})}

\providecommand{\CyclicGroup}[1]{\mathcal C_{#1}}
\providecommand{\ElementaryAbelianGroup}[2]{\CyclicGroup{#1}^{#2}}%{\mathcal E_{#1}({#2})} % first param: length of cycle, second param: number of cycles, i.e. we have \CyclicGroup{#1} \times \dots \times \CyclicGroup{#1} with #2 factors.
\providecommand{\DihedralGroup}[1]{\mathcal D_{#1}}
\providecommand{\AlternatingGroup}[1]{\mathcal A_{#1}}
\providecommand{\SymmetricGroup}[1]{\mathcal S_{#1}}
\providecommand{\SemiDirectProductGroup}[3]{\ElementaryAbelianGroup{#1}{#2} \rtimes \CyclicGroup{#3}} % first param: characteristic, second param: Number of C_#1-s, third param: order of additional cyclic group
\providecommand{\PSL}[2]{\mathrm{PSL}_{#1}({#2})} % first param: matrix dimension, second param: field or order of the field.
\providecommand{\PGL}[2]{\mathrm{PGL}_{#1}({#2})} % first param: matrix dimension, second param: field or order of the field.

\providecommand{\typeOfField}[2]{\mathrm{F}[{#1},{#2}]}
\providecommand{\typeOfFieldAutSubGroup}[2][]{\mathrm{U}_{#1}({#2})}

\providecommand{\GF}[1]{\mathbb F_{#1}}

\providecommand{\FixedField}[2]{{#1}^{{#2}}} % #1: Körper, #2: Automorphismengruppe
\providecommand{\AlgebraicClosure}[1]{\overline {#1}} % algebraischer Abschluss von #1

\providecommand{\Jacobian}[1][]{{\mathbb J}_{#1}}

\newcounter{algorithmStepCount}

\newcounter{algorithmSubStepCount}
               {\setcounter{algorithmStepCount}{1}%
	        \par\mbox{}%
	        \begin{list}%
		   {(\arabic{algorithmStepCount}.\arabic{algorithmSubStepCount})}%
 	           { \usecounter{algorithmSubStepCount}%
                   }%
	       }%
	       {\end{list}\mbox{}\par}

\providecommand{\valuation}[1]{\mathrm {v}_{#1}} % arg: place
\providecommand{\principalDivisor}[2][]{({#2})^{#1}} % first (optional) arg=function field. second arg=element.
\providecommand{\zeroDivisor}[2][]{({#2})^{#1}_{0}}
\providecommand{\poleDivisor}[2][]{({#2})^{#1}_{\infty}}

\providecommand{\AdeleSpace}[2][\empty]{\ifx#1\empty{\mathcal {A}_{#2}}\else{\mathcal {A}_{#2}(#1)}\fi}
\providecommand{\WeilDifferentials}[2][\empty]{\ifx#1\empty{{\Omega}_{#2}}\else{{\Omega}_{#2}(#1)}\fi}

\providecommand{\support}{\mathrm {supp}} % support of some divisor.

\newtheorem{theorem}{Theorem}
\newtheorem{proposition}[theorem]{Proposition}
\newtheorem{lemma}[theorem]{Lemma}
\newtheorem{corollary}[theorem]{Corollary}

\theoremstyle{definition}
\newtheorem{example}{Example}
\newtheorem{remark}[example]{Remark}
\newtheorem{definition}[example]{Definition}
\newtheorem{algorithm}[example]{Algorithm}

\begin{document}
\title[Automorphism Groups of Hyperelliptic Function Fields]
      {Computing the Automorphism Groups of Hyperelliptic Function Fields}
\author{Norbert G\"ob}
\thanks{The author is partially supported by  Fraunhofer~ITWM, Kaiserslautern and
               BGS~AG, Mainz.}
\date{May 20, 2003}
\maketitle

\section{Introduction}
The purpose  of this paper  is to propose an  efficient method to  compute the
automorphism  group  of  an  arbitrary hyperelliptic  function  field  over  a
given  ground field  of  characteristic $>2$  as well  as  over its  algebraic
extensions. Beside theoretical applications, knowing the automorphism group of
a hyperelliptic function field also is useful in cryptography:

The  Jacobians  of  hyperelliptic  curves   have  been  suggested  by  Koblitz
as   groups   for  cryptographic   purposes,   because   the  computation   of
the   discrete   logarithm   is   believed   to   be   hard   in   this   kind
of   groups    (\cite{koblitz:1989}).   In   order   to    obtain   ``secure''
Jacobians   it  is   necessary  to   prevent  attacks   like  Pohlig/Hellman's
(\cite{pohligHellman:1978}),    Frey/R\"uck's   (\cite{freyRueck:1994})    and
Duursma/Gaudry/Morain's  (\cite{duursmaGaudryMorain:1999}). The  latter attack
is only feasible,  if the corresponding function field has  an automorphism of
large order. To forestall the Pohlig-Hellman  attack, one needs to assert that
the group  order is  almost prime,  i.e.\ it  ought to  contain a  large prime
factor  $p_0$. To  prevent  the  Frey-R\"uck attack,  $p_0$  needs to  possess
additional properties.

Therefore, one needs to know both the automorphism group of the function field
and the order of the Jacobian.  Unfortunately, there is no efficient algorithm
known  to compute  this order  for  arbitrary hyperelliptic  curves. Only  for
specific types of curves, divisor class counting\footnote{i.e.\ computation of
the order  of the Jacobian}  is feasible for cryptographically  relevant group
sizes (e.g.\ \cite{sakaiSakuraiIshizuka:1998}, \cite{gaudryHarley:2000}).

A  theorem  by  Madan  (\cite{madan:1970})  implies  that  $\cardinality{
\Jacobian[F]}$  divides $\cardinality{  \Jacobian[F']}$ whenever  $F \subseteq
F'$ is  a (hyper-)elliptic subfield  of a hyperelliptic function  field s.th.\
\mbox{$[F':F] <  \infty$}. Thus,  a hyperelliptic  function field  with secure
Jacobian  will  most likely  have  a  trivial  automorphism group,  i.e.\  one
consisting  of the  hyperelliptic involution,  only. Therefore,  the proposed technique
provides a quick test to check whether a given hyperelliptic curve may yield a
secure Jacobian,  i.e.\ whether  it is worthwhile  to apply  expensive divisor
class counting algorithms.

Let us  outline the afore  mentioned algorithm briefly.  It is well  known that
the  automorphism group  of a  hyperelliptic  function field  is finite  (cf.\
\cite{schmid:1938}). For  each finite  group, which can  occur as  subgroup of
such an  automorphism group, Brandt gave  a normal form for  the corresponding
hyperelliptic function  fields and  explicit formulas for  these automorphisms
(cf.\ \cite{brandt:1988}). Brandt's results only apply to function fields over
algebraically closed constant fields, but this  is no hindrance as we will see
later. For now, we suppose the constant field to be algebraically closed.

Hence, computing  the automorphism  group reduces to  the problem  of deciding,
whether  a  given  hyperelliptic  function   field  has  a  defining  equation
of  the  form   given  by  Brandt's  theorems.  This  can   be  checked  using
theorem~\ref{thm:relationOfGeneratorsOfHyperellipticFields}, which states that
two hyperelliptic  function fields $k(t,u)$,  $k(x,y)$ with $u^2=D_t$, $y^2=D_x$
are equal  iff $x=\frac
{\alpha_0t +\alpha_1} {\alpha_2t + \alpha_3}$  for some $\alpha_i \in {k}$ and
$y=\varphi  u$,  where  $\varphi  \in  {k}(t)$ can  be  determined  from  the
$\alpha_i$.  Hence, we  substitute $x=\frac  {\alpha_0t +\alpha_1}  {\alpha_2t
+\alpha_3}$  symbolically into  $D_x$.  Computing $\varphi$  according to  the
theorem  and  comparing  coefficients  of $D_t$  and  $\varphi^{-2}D_x  (\frac
{\alpha_0t +\alpha_1} {\alpha_2t  +\alpha_3}) = \varphi^{-2} y^2  = u^2 =D_t$,
we obtain  a system of polynomial  equations for the $\alpha_i$.  These can be
tested for solvability or even solved using Gr\"obner basis methods.

If     the     constant     field     $k$     is     algebraically     closed,
algorithm~\ref{alg:compAutWrtAlgClosedK}  seems  to   be  the  only  efficient
possibility  known  to   compute  the  automorphism  group   of  an  arbitrary
hyperelliptic  function  field.  For  finite  $k$,  the  method  described  in
section~\ref{subsec:compAut:arbitraryConstants} is an  alternative approach to
the \texttt{AutomorphismGroup} function in \cite{stoll:2001}.

\section{Notation and Fundamental Facts}

Throughout this paper, we use  the notations from \cite{stichtenoth:1993}. For
the  reader's convenience,  we  recall the  essential  notations: The  natural
numbers $\N$ start at $0$, $\Nplus:=  \N \setminus \{0\}$. The greatest common
divisor of two  integers or polynomials $p,q$ is denoted  by $\gcd{p}{q}$. The
unit  group  of a  field  $k$  is denoted  by  $\unitGroup{k}  := k  \setminus
\{0\}$.  Let  $k$ be  some  field  of  characteristic  $p>2$, and  $g\in  \N$,
$g>1$.  A  {\em  hyperelliptic  function  field}\/ of  genus  $g$  over  $k$  is
defined  to be  a  field $F:=k(x,y)$  s.th.\ $x$  is  transcendental over  $k$
and  $y^2=D(x)$,  where $D  \in  k[x]$  is  a  monic separable  polynomial  of
degree $2g+1$  or $2g+2$.  The {\em  automorphism group of  $F$}\/ is  the group
$\Automorphismgroup  {F}  {k}$  of  field automorphisms  of  $F$  fixing  $k$.
If  $U  \le  \Automorphismgroup  {F}  {k}$,  we  denote  the  fixed  field  of
$U$  by $\FixedField  {F}{U}$.  The algebraic  closure of  $k$  is denoted  by
$\AlgebraicClosure{k}$. If $P$ is a  place of $F$, $\valuation{P}$ denotes the
valuation corresponding to $P$. For $t  \in F$ we denote the principal divisor
of $t$  by $\principalDivisor{t}$, its  zero divisor by  $\zeroDivisor{t}$ and
its pole  divisor by  $\poleDivisor{t}$. If $\poleDivisor{t}$  is a  place, we
also denote it by $\infty_t :=  \poleDivisor{t}$ and call it the {\em infinite
place w.r.t.\ $t$}.

Our  aim is  to  compute the  automorphism group  of  any given  hyperelliptic
function field  $k(x,y)$, $y^2=D$. As  mentioned above, Brandt  gives normal
forms of  hyperelliptic function fields  for each possible finite  subgroup of
the  automorphism group  (cf.\  Brandt's  Ph.D.\ thesis,  \cite{brandt:1988}).
Since the automorphism  group of such a field is a central extension
of  $\AutomorphismGroup{k(x)}{k}$ by  the $\CyclicGroup{2}$  generated by  the
hyperelliptic involution, Brandt rather investigates the possible subgroups of
$\AutomorphismGroup{k(x,y)}{k} / \CyclicGroup{2}$,  i.e.\ he characterizes the
fields by their ``type'' which is defined as follows.

\begin{definition}{Type of field}{$\typeOfField{G}{k}$}
  Let $F/k$ be a hyperelliptic function field and $G$ some finite group. $F$
  is called  a function field  of {\em type $\typeOfField{G}{k}$},  if there
  are  finite groups  $C, U$,  s.th.\ $U  \le \Automorphismgroup{F}{k}$,  $C
  \unlhd U$, $C \cong \CyclicGroup{2}$, $\FixedField{F}{C}$ is a rational function field
  over $k$, and $U/C \cong G$.

  We         denote         such         a        group         $U$         by
  $\typeOfFieldAutSubGroup     {G}$\index{\typeOfFieldAutSubGroup[F]{G}}    or
  $\typeOfFieldAutSubGroup[F]{G}$,  although  $U$  needs not  to  be  uniquely
  determined by $F$, $G$ and $k$. We will only use this notation to state that
  a specific group can be used as $U$ in this definition.
  % Problem: U ist nicht notwendig eindeutig bestimmt!

  For  extension  fields  $k'  \supseteq  k$,  we  call  $F$  to  be  of  type
  $\typeOfField {G}{k'}$ iff the constant  field extension $Fk'/k'$ is of type
  $\typeOfField {G}{k'}$.
\end{definition}

The   following   types   can   occur  for   hyperelliptic   function   fields
over    algebraically    closed     constant    fields    of    characteristic
$p>2$:   $\typeOfField{\CyclicGroup{n}}{k}$,    where   $\gcd{n}{p}    =   1$,
$\typeOfField{\ElementaryAbelianGroup{p}{m}}{k}$   for    some   $m\in\Nplus$,
$\typeOfField  {\DihedralGroup {n}}  {k}$, where  $\gcd{n}{p} =  1$ or  $n=p$,
$\typeOfField {\AlternatingGroup  {4}}{k}$, $\typeOfField  { \AlternatingGroup
{5}  }{k}$,   $\typeOfField  {  \SymmetricGroup  {4}   }  {k}$,  $\typeOfField
{  \SemiDirectProductGroup{p}{m}{n}}{k}$,  where  $\gcd{n}{p}=1$  and  $m  \in
\Nplus$,  $\typeOfField   {  \PSL{2}{p^m}   }  {k}$,  where   $m\in\Nplus$  and
$\typeOfField  {  \PGL{2}{p^m} }  {k}$,  where  $m\in  \Nplus$. As  one  needs
to  consider several  cases  for  some of  these  types,  the theorem  stating
Brandt's  normal  forms  contains  14   case  distinctions.  For  brevity,  we
only  consider the  types with  the  smallest and  largest possible  subgroups
as  well  as  a  subgroup  which  we will  need  in  our  examples.  Hence, we
restrict  ourselves  to  the types  $\typeOfField{\CyclicGroup{n}}{k}$,  where
$\gcd{n}{p}  = 1$,  $\typeOfField{\ElementaryAbelianGroup{p}{m}}{k}$ for  some
$m\in\Nplus$,  $\typeOfField {\DihedralGroup{n}}  {k}$, where  $\gcd{n}{p}=1$,
and $\typeOfField  { \PGL{2}{p^m} }  {k}$, where $m\in \Nplus$.  The remaining
cases are  similar to these and  can be found in~\cite{goeb:2003}  or directly
in~\cite{brandt:1988}.

\begin{theorem}[Brandt]\label{thm:brandtsNormalForms}\index{Brandt's normal forms}
  Let $F$  be a  hyperelliptic function field  over an  algebraically closed
  constant field $k$ of characteristic $p\ge 3$. 
  Then the types of $F$ are characterized as follows
  \begin{enumerate}
  \item
    $F$  is of  type $\typeOfField{\CyclicGroup{n}}{k}$  for $n\in  \Nplus$ with
    $\gcd{n}{p}=1$ iff there are $t,u \in F$, s.th.\ $F=k(t,u)$, $u^2= t^\nu
    \prod_{j=1}^s(x^n-a_j)$, where $\nu \in \{0,1\}$, $s\in \N$ and the $a_j
    \in \unitGroup{k}$ are pairwise distinct.

    In this case, $\typeOfFieldAutSubGroup[F]{\CyclicGroup{n}}$ is generated
    by $\varphi: t \mapsto t$, $u \mapsto -u$ and $\psi: t\mapsto \eta^2 t$,
    $u \mapsto  \eta^\nu u$,  where $\eta$  is a  primitive $2n$-th  root of
    unity.
  \item
    $F$  is  of type  $\typeOfField{\ElementaryAbelianGroup{p}{m}}{k}$  with
    $m\in \Nplus$ iff there are $t,u  \in F$ and a subgroup $A$ of the
    additive group of $k$ of order $\cardinality{A}=p^m$, s.th.\ $F=k(t,u)$,
    $u^2=\prod_{j=1}^s\left(\prod_{a\in A} (x+a)  - a_j\right)$, where $s\in
    \N$ and the $a_j \in k$ are pairwise distinct.

    In this case,  $\typeOfFieldAutSubGroup [F] {\ElementaryAbelianGroup {p}
    {m}}$  is generated  by  $\varphi:  t \mapsto  t$,  $u  \mapsto -u$  and
    all $\psi_a: t\mapsto t+a$, $u \mapsto u$ with $a \in A$.
  \item
    $F$ is  of type  $\typeOfField {\DihedralGroup {n}}  {k}$, where  $n \in
    \Nplus$, $\gcd  {n} {p} =1$ iff there  are $t,u \in
    F$, s.th.\ $F=k(t,u)$, $u^2  = t^{\nu_0} (t^n-1)^{\nu_1} (t^n+1)^{\nu_2}
    \prod_{j=1}^s (t^{2n} - a_j t^n +  1)$, where $\nu_j \in \{0,1\}$, $s\in
    \N$ and the $a_j \in  k \setminus \{ \pm 2  \}$ are pairwise
    distinct. If $n=2$ or $n \equiv 1 \mod 2$, we need to have $\nu_1=\nu_2$.

    In  this case,  $\typeOfFieldAutSubGroup  [F]  {\DihedralGroup {n}}$  is
    generated by  $\varphi: t \mapsto t$,  $u \mapsto -u$, $\psi:  t \mapsto
    \eta^2 t$, $u \mapsto \eta^{\nu_0}u$ and $\sigma: t \mapsto \frac 1t$, $u
    \mapsto \frac {i^{\nu_1}u}  {t^m}$, where $\eta$ is  a primitive $2n$-th
    root of unity, $i^2=-1$ and $m = \frac  12 n (\nu_1 + \nu_2) + 2 \nu_0 +
    2ns$.
  \item[(14)]
    $F$  is of  type  $\typeOfField  { \PGL{2}{p^m}  }  {k}$  iff there  are
    $t,u  \in F$,  s.th.\ $F=k(t,u)$,  
    \begin{multline*}
      u^2  = (t^r -  t)^{\nu_0} 
	     \left(  (t^r -  t)^{r-1} +  1 \right)^{\nu_1} 
	     \\
	     \prod_{j=1}^s  \left( \left(  (t^r - t)^{r-1}  +  1  \right)^{r+1} 
		-  a_j  (t^r  -  t)^{r^2-r}  \right),
    \end{multline*}
    where  $\nu_j \in  \{0,  1\}$, $s  \in  \N$, $r=p^m$  and  the $a_j  \in
    \unitGroup{k}$ are pairwise distinct.

    In this case, $\typeOfFieldAutSubGroup  [F] {\PGL{2}{p^m}}$ is generated
    by $\varphi: t  \mapsto t$, $u \mapsto -u$, $\psi:  t \mapsto \eta^2 t$,
    $u \mapsto \eta^{\nu_0}  u$, $\sigma: t \mapsto t+1$, $u  \mapsto u$ and
    $\tau: t \mapsto \frac 1t$, $u  \mapsto \frac u{t^n}$, where $\eta$ is a
    primitive $2(p^m -1)$-th root of unity and 
    \[
      n = \frac 12 \left(  (p^m +1) \nu_0  + p^m (p^m -1) \nu_1 
	  + p^m(p^{2m} - 1)s \right).
    \]
  \end{enumerate}
\end{theorem}
\begin{proof}
  A  slightly more  general theorem  is proved  by Rolf  Brandt in  his Ph.D.\
  thesis \cite{brandt:1988}:  He characterizes the types  of cyclic extensions
  of rational  function fields over  algebraically closed constant  fields. We
  list the references  for each of the  stated facts, citing the  proof that a
  function field of the given type has the given normal form, first. The proof
  of the inverse implication and the generators are given thereafter.
  \begin{enumerate}
  \item % \CyclicGroup {n}.
    \cite[Satz~5.1]{brandt:1988},
    \cite[Satz~5.6]{brandt:1988} and
    \cite[Lemma~5.5]{brandt:1988}.
  \item % \ElementaryAbelianGroup {p} {m}.
    \cite[Satz~6.3]{brandt:1988} and its proof.
  \item % \DihedralGroup{n}, (n,p)=1, n==0 mod 2.
    Cf.\ \cite[Satz~7.3]{brandt:1988},
    \cite[Satz~7.5]{brandt:1988} and
    \cite[Lemma~7.4]{brandt:1988},
    in the case $n \equiv 0 \mod 2$.
    Otherwise, we apply \cite[Satz~7.9]{brandt:1988}, as $p\ge
    3$  and  $\gcd{n}{p} =1$  obviously  imply  $\gcd  {2n}  {p} =  1$.  The
    generators  and  the  inverse  implication  are  proved  analogously  to
    \cite[Satz~7.5]{brandt:1988} and \cite[Lemma~7.4]{brandt:1988}.
  \item[(14)] % \PGL{2}{p^m}
    \cite[Satz~13.1]{brandt:1988},
    \cite[Satz~13.6]{brandt:1988} and
    \cite[Lemma~13.2]{brandt:1988}.
  \end{enumerate}
\end{proof}

Let us illustrate this theorem and the related problems with an example.
\begin{example}\label{ex:brandtsNormalFormPartiallyVisible}
  We   consider  $F:=\AlgebraicClosure{\GF{7}}(x,y)$,
  \[
     y^2=   x^5+x^3+x  
      = x(x+2)(x-2)(x+3)(x-3) = x(x^2-4)(x^2-2).
  \]
  Obviously    $F$     is    of    type     $\typeOfField    {\CyclicGroup{2}}
  {\AlgebraicClosure{\GF{7}}}$.  The  basis  $x,y$  of  $F$  is  not  uniquely
  determined by  $F$, neither  is the defining  equation. Therefore,  we cannot
  immediately see if $F$ is of any other types.
\end{example}

In the  following section, we solve  this problem, i.e.\ we propose an efficient
possibility to find  out, if a hyperelliptic function  field has a given
normal form.

\section{Relations Between Bases}\label{sec:defEq}

In      this     section      we     show      the     connection      between
different      bases      of      a     hyperelliptic      function      field
(cf.\     theorem~\ref{thm:relationOfGeneratorsOfHyperellipticFields}):     If
$k(t,u)=k(x,y)$  is  a  hyperelliptic  function  field,  then  $x$  needs  to
be   a   fraction   of   linear   polynomials  in   $t$   and   the   relation
between   $u$   and   $y$   can   be   computed   easily  from these polynomials.
In  contrast   to  theorem~\ref{thm:brandtsNormalForms}, we
do   not   need   to   have    an   algebraically   closed   constant   field, here;
theorem~\ref{thm:relationOfGeneratorsOfHyperellipticFields}     applies     to
hyperelliptic function fields over arbitrary constant fields of characteristic
$\neq 2$.
This  theorem is  one   of  the   core  components   of  our  algorithm  for
computing  the  automorphism group  of  a  hyperelliptic function  field, as
we will see in section~\ref{sec:compAut}.

\subsection{Relations Between the Variable Symbols}

Here, we show that $x$ can be  represented as a fraction of linear polynomials
in $t$. We start our proof by citing the following lemma:

\begin{lemma}\label{lem:kxEqKtInHyperellFields}
  Let  $k(t,u)=k(x,y)$   be  a  hyperelliptic  function   field,  $u^2=D_t$,
  $y^2=D_x$, where  $D_t \in k[t]$  and $D_x  \in k[x]$ are  separable monic
  polynomials. Then $k(t)=k(x)$.
\end{lemma}
\begin{proof}
  \cite[Proposition~VI.2.4]{stichtenoth:1993}.
\end{proof}
% TODO: Gilt diese Aussage auch f\"ur g=1?

Lemma~\ref{lem:kxEqKtInHyperellFields}    means,    that   the    following
proposition can be applied to our situation, i.e.\ in hyperelliptic function
fields with two given bases, we always  have $k(x)=k(t)$. We see that $x$ is
a fraction of linear polynomials in $t$ in this case:

\begin{proposition}\label{prop:generatorsOfRationalFuncFields}
  Let   $k(t)$   be  a   rational   function   field   and  $x   \in   k(t)$
  s.th.\   $k(t)=k(x)$.   Then   there  are   $\alpha_0,   \dots,   \alpha_3
  \in   k$    with   $x=\frac{\alpha_0t+\alpha_1}{\alpha_2t+\alpha_3}$   and
  $\alpha_0\alpha_3-\alpha_1\alpha_2 \neq 0$.
\end{proposition}
\begin{proof}
  As $x \in k(t)$, there are polynomials $\varphi, \psi \in k[t]$, s.th.\ $x
  = \frac{\varphi}{\psi}$  and $\gcd{\varphi}{\psi} \in k$.  We consider the
  principal divisor of $x$. \cite[Theorem I.4.11]{stichtenoth:1993} implies
  \[
    \deg(\zeroDivisor{x}) = \deg(\poleDivisor{x}) =[k(t):k(x)] =1.
  \]
  Let us consider the  case $\infty_t \notin \support{\principalDivisor{x}}$,
  first.   Then   $0=\valuation{\infty_t}(x)=\deg_t(\psi) - \deg_t(\varphi)$,
  i.e.\   $\deg_t(\varphi)   =
  \deg_t(\psi)$. As $\varphi, \psi \in k[t]$, we get $\poleDivisor {\varphi}
  =  \deg_t  (\varphi) \infty_t  =  \deg_t  (\psi) \infty_t  =  \poleDivisor
  {\psi}$.  We  have  $\principalDivisor{x} =  \principalDivisor{\varphi}  -
  \principalDivisor{\psi} = \zeroDivisor{\varphi}  - \poleDivisor{\varphi} -
  (  \zeroDivisor{\psi} -  \poleDivisor{\psi}  )  = \zeroDivisor{\varphi}  -
  \zeroDivisor{\psi}$, i.e.\  $\zeroDivisor{x} =  \zeroDivisor{\varphi}$ and
  $\poleDivisor{x} = \zeroDivisor{\psi}$. Thus,
  \begin{multline*}
    \deg_t(\varphi)   =\deg(\poleDivisor{\varphi})
    =\deg(\zeroDivisor{\varphi})
    =\deg(\zeroDivisor{x}) \\
    =1      =\deg(\poleDivisor{x})
    =\deg(\zeroDivisor{\psi}) = \deg(\poleDivisor{\psi}) =\deg_t(\psi).
  \end{multline*}
  Thus  there  are  $\alpha_i \in  k$  s.th.\  $\varphi=\alpha_0t+\alpha_1$,
  $\psi=\alpha_2t+\alpha_3$ and $\alpha_0\alpha_3 - \alpha_1\alpha_2 \neq 0$
  as claimed.

  If   $\infty_t    \in   \support{\principalDivisor{x}}$,    we   obviously
  have   $\deg_t(\varphi)   \neq    \deg_t(\psi)$.   W.l.o.g.\   we   assume
  $\valuation{{\infty_t}}(x)<0$     (consider    $\frac{1}{x}$     in    the
  other   case).    As   $\deg(\poleDivisor{x})=1$,   we   need    to   have
  $\valuation{{\infty_t}}(x)=-1$. Thus
  \[
     -1  =  \valuation {{\infty_t}}  (x) = \deg_t (\psi) - \deg_t (\varphi),
  \]
  i.e.\   $\deg_t(\psi)=\deg_t(\varphi)-1$.   As  $\poleDivisor{\varphi}   =
  \deg_t(\varphi)   \infty_t$   and   $\poleDivisor{\psi}   =   \deg_t(\psi)
  \infty_t$, we infer
  \[
    \principalDivisor{x}        =        \principalDivisor{\varphi}        -
    \principalDivisor{\psi} =  \zeroDivisor{\varphi} - \poleDivisor{\varphi}
    -  \zeroDivisor{\psi}  +  \poleDivisor{\psi} =  \zeroDivisor{\varphi}  -
    \zeroDivisor{\psi} - \infty_t.
  \]
  Thus,        we       have        $\zeroDivisor{x}=\zeroDivisor{\varphi}$,
  i.e.\      
  \[
    \deg_t(\varphi) =  \deg(\poleDivisor{\varphi}) = \deg(\zeroDivisor{\varphi})
    = \deg(\zeroDivisor{x}) = 1.
  \]
  Furthermore $\deg_t(\psi)=\deg_t(\varphi)-1 = 0$.
  We   obtain   $x   =   \frac  {\varphi}   {\psi}   =   \frac{\alpha_0t   +
  \alpha_1}{\alpha_3}$ with  $\alpha_i \in k$ and $\alpha_0\alpha_3  \neq 0$ as
  claimed.
\end{proof}

Summing up these facts, we obtain, that $x$ is a fraction of linear polynomials
in~$t$, if $k(x,y)=k(t,u)$:
\begin{corollary}\label{cor:relationOftAndxinKtuAndKxy}
  % TODO: allow g=1 if lem:kxEqKtInHyperellFields is true in this case
  Let  $k(t,u)=k(x,y)$   be  a  hyperelliptic  function   field,  $u^2=D_t$,
  \mbox{$y^2=D_x$},  where  $D_t  \in  k[t]$  and  $D_x  \in  k[x]$  are  separable
  monic   polynomials.   Then   there   are   $\alpha_0,   \dots,   \alpha_3
  \in   k$    with   $x=\frac{\alpha_0t+\alpha_1}{\alpha_2t+\alpha_3}$   and
  $\alpha_0\alpha_3-\alpha_1\alpha_2 \neq 0$.
\end{corollary}
\begin{proof}
  By lemma~\ref{lem:kxEqKtInHyperellFields}, we have $k(t)=k(x)$. Thus
  proposition~\ref{prop:generatorsOfRationalFuncFields} implies the
  existence of the $\alpha_i$.
\end{proof}

\subsection{Relation Between the Square Roots}

Since we know,  how $t$ and $x$ are related  in a hyperelliptic function
field for which  we have two bases $k(t,u)=k(x,y)$, we  proceed studying the
relationship between  $u$ and $y$.  The next lemma tells  us, that $y$  is a
multiple of $u$ over $k(t)$:

\begin{lemma}\label{lem:relationOfuAndyinKtuAndKxyLinear}
  Let $F=k(t,u)=k(x,y)$,  $u^2=D_t$, $y^2=D_x$  be a  hyperelliptic function
  field over  a constant field  $k$ of  characteristic $\neq 2$,  where both
  $D_t \in  k[t]$ and $D_x \in  k[x]$ are monic separable  polynomials. Then
  there is some $\varphi \in k(t)\setminus \{0\}$, s.th.\ $y= \varphi u$.
\end{lemma}
\begin{proof}
  As  $y \in  F=k(t,u)$ and  $[k(t,u):k(t)]=2$, there  are $\varphi,  \psi \in
  k(t)$ s.th.\  $y = \varphi  u + \psi$. Let  us suppose $\varphi=0$.  Then we
  had  $y  \in  k(t)$.  From  lemma~\ref{lem:kxEqKtInHyperellFields}  we  know
  that  $k(x)=k(t)$. Thus  we had  $y\in k(x)$,  i.e.\ $k(x,y)=k(x)$  implying
  $[k(x,y):k(x)]=1$, which contradicts  $[k(x,y):k(x)]=2$. Therefore  $\varphi
  \neq 0$.

  Substituting our representation of $y$ into its minimal polynomial we get
  \[
     D_x  =  y^2  =  (\varphi  u +  \psi)^2  =  \varphi^2u^2  +
     2\varphi\psi u + \psi^2 =  \varphi^2D_t + 2\varphi\psi u +
     \psi^2.
  \]
  Thus  $2\varphi\psi  u  \in  k(t)=k(x)$.  As $u  \notin  k(t)$,  this  leads
  to  $2  \varphi  \psi  =  0$,  from  with  we  conclude  $\psi  =0$  because
  $\characteristic{k}\neq 2$ and $\varphi \neq 0$.
\end{proof}

Knowing that $y=\varphi  u$, we will examine $\varphi$ more  closely. We start
by the following lemma,  which is quite technical, but will  be useful in the
subsequent proofs.

\begin{lemma}\label{lem:relationOfuAndyinKtuAndKxyFactorsOfDxWrtt}
  Let $k(t,u)=k(x,y)$,  $u^2=D_t$, $y^2=D_x$  be a  hyperelliptic function
  field
  over a constant field $k$ of  characteristic $\neq 2$, where both $D_t \in
  k[t]$ and $D_x \in k[x]$ are monic separable polynomials.
  % CAUTION: Separability is needed, here. If a square free polynomial
  % is not necessarily separable, we need to suppose separability!
  Let      $x=\frac{\alpha_0t+\alpha_1}{\alpha_2t+\alpha_3}$,      $\alpha_i
  \in     k$,     $\alpha_0\alpha_3-\alpha_1\alpha_2      \neq     0$     as
  stated      in     corollary~\ref{cor:relationOftAndxinKtuAndKxy}      and
  $y=\varphi     u$,     $\varphi     \in    k(t)\setminus     \{0\}$     as
  in     lemma~\ref{lem:relationOfuAndyinKtuAndKxyLinear}.    Then     there
  are    $d_x:=\deg_x(D_x)$    pairwise    relatively   prime    $p_i    \in
  \AlgebraicClosure{k}[t]$, $\deg_t(p_i)\le 1$ s.th.\
  \[
    D_t = \varphi^{-2} (\alpha_2t+\alpha_3)^{-d_x}
         \prod_{i=1}^{d_x} p_i.
  \]
  Furthermore we have
  \begin{enumerate}
  \item
    $p_i = (\alpha_0 + \alpha_2\eta_i)t  + \alpha_1 - \alpha_3\eta_i$, where
    $\eta_i \in \AlgebraicClosure{k}$,  $i=1, \dots, d_x$ are  the zeroes of
    $D_x$.
  \item
    $d_x-1 \le \deg_t\left( \prod_{i=1}^{d_x} p_i \right) \le d_x$.
  \item
    Let  $q  \in  \AlgebraicClosure{k}[t]$  be  linear.  Then  $q^2  \ndivides
    \prod_{i=1}^{d_x} p_i$. In particular, $(\alpha_2t + \alpha_3)^2 \ndivides
    \prod_{i=1}^{d_x} p_i$.
  \end{enumerate}
\end{lemma}
\begin{proof}
  We    factor   $D_x$    over   $\AlgebraicClosure{k}$    into
  $D_x   =    \prod_{i=1}^{d_x}   (x-\eta_i)$,    $\eta_i   \in
  \AlgebraicClosure{k}$,  $\eta_i \neq  \eta_j$ for  all $i\neq
  j$. This yields
  \begin{align*}
    D_t=&u^2
            = \varphi^{-2} y^2
            = \varphi^{-2} D_x
            = \varphi^{-2} \prod_{i=1}^{d_x} (x-\eta_i)
         \\
            =& \varphi^{-2} \prod_{i=1}^{d_x} (\frac{\alpha_0t+\alpha_1}
                                             {\alpha_2t+\alpha_3}-\eta_i)
         \\
       =& \varphi^{-2}
       \prod_{i=1}^{d_x} \frac{(\alpha_0-\alpha_2\eta_i)t
                               +\alpha_1-\alpha_3\eta_i}
                       {\alpha_2t+\alpha_3}
         \\
       =:& \varphi^{-2}
      \prod_{i=1}^{d_x} \frac{p_i}{\alpha_2t+\alpha_3}
         \\
       =& \varphi^{-2} (\alpha_2t+\alpha_3)^{-d_x}
       \prod_{i=1}^{d_x} {p_i} \in k[t],
  \end{align*}
  where  $p_i  := (\alpha_0  -  \alpha_2\eta_i)t  + \alpha_1  -
  \alpha_3\eta_i$. Suppose there were indices $i \neq j$ s.th.\
  $p_i$  and  $p_j$ had  a  common  divisor of  nonzero  degree
  w.r.t.  $t$. Then  we  had $p_i=\beta  p_j$  for some  $\beta
  \in \AlgebraicClosure{k} \setminus  \{0\}$, i.e.\ $(x-\eta_i)
  =   \frac{p_i}   {\alpha_2t+\alpha_3}  =   \beta   \frac{p_j}
  {\alpha_2t+\alpha_3} = \beta(x-\eta_j)$. Thus, $D_x$ were not
  separable. Contradiction.  Therefore, the $p_i$  are pairwise
  relatively prime, which proves our main claim.

  Let us proceed by examining the supplementary statements.
  Obviously,
  \[ 
    \deg_t \left (\prod_{i=1}^{d_x}  p_i  \right)   \le  d_x.
  \]
  If  $\deg_t\left(\prod_{i=1}^{d_x} p_i  \right) <  d_x-1$, there  were two
  indices  $i\neq  j$  s.th.\  $p_i,  p_j  \in  \AlgebraicClosure{k}$,  thus
  $\alpha_0  -  \alpha_2   \eta_i  =  \alpha_0  -  \alpha_2   \eta_j  =  0$,
  i.e.\  $\alpha_0=\alpha_2 \eta_i  = \alpha_2  \eta_j$. Hence,  $\alpha_2(\eta_i
  -\eta_j)=0$ which yields  $\alpha_2=0$ since $\eta_i \neq  \eta_j$. Now we
  can easily deduce $\alpha_0=0$ from  $\alpha_0 - \alpha_2\eta_i =0$. Since
  $\alpha_0\alpha_3-\alpha_1\alpha_2  \neq 0$,  this is  not possible.  Thus
  $\deg_t\left(\prod_{i=1}^{d_x} p_i \right) \ge d_x-1$.

  Finally, let $q^\nu \divides \prod_{i=1}^{d_x} p_i$ for some linear
  $q\in k[t]$  and   $\nu  \in \Nplus$. As  $\deg_t(p_i)\le 1$,
  there are $\nu$  factors $p_{i_1}, \dots, p_{i_\nu}$,  which are multiples
  of  $q$. Thus  $p_{i_1}, \dots,  p_{i_\nu}$ are  scalar multiples  of each
  other. If $\nu>1$,  this contradicts the relative primality  of the $p_i$.
  This proves the last claim.
\end{proof}

The following lemma states, that $\varphi^{-1}$ is a non-zero polynomial in $t$:

\begin{lemma}\label{lem:relationOfuAndyinKtuAndKxyLinFactorIsPolyInv}
  Let  $k(t,u)=k(x,y)$,  $u^2=D_t$,  $y^2=D_x$  be  a  hyperelliptic  function
  field  over a  constant field  $k$ of  characteristic $\neq  2$, where  both
  $D_t  \in  k[t]$  and  $D_x  \in  k[x]$  are  monic  separable  polynomials.
  Let  $x=\frac  {\alpha_0t +  \alpha_1}  {\alpha_2t  + \alpha_3}$,  $\alpha_i
  \in  k$,  $\alpha_0\alpha_3   -  \alpha_1\alpha_2  \neq  0$   as  stated  in
  corollary~\ref{cor:relationOftAndxinKtuAndKxy}  and  $y=\varphi   u$  as  in
  lemma~\ref{lem:relationOfuAndyinKtuAndKxyLinear}. Then we have $\varphi^{-1}
  \in k[t] \setminus \{0\}$.
\end{lemma}
\begin{proof}
  Lemma~\ref{lem:relationOfuAndyinKtuAndKxyFactorsOfDxWrtt}  implies $D_t  =
  \varphi^{-2} (\alpha_2t+\alpha_3)^{-d_x} \prod_{i=1}^{d_x} {p_i}$. Suppose
  $\varphi^{-1} =  \frac {\varphi_1} {\varphi_0}  \notin k[t]$. As  $D_t \in
  k[t]$, $\varphi_0^2$  needs to be canceled  by $\prod_{i=1}^{d_x} {p_i}$.
  Let $q  \in \AlgebraicClosure{k}[t]$  be a  linear factor  of $\varphi_0$.
  Thus   $q^2   \divides   \prod_{i=1}^{d_x}   {p_i}$,   which   contradicts
  lemma~\ref{lem:relationOfuAndyinKtuAndKxyFactorsOfDxWrtt}. Thus we need to
  have $\varphi^{-1} \in k[t]$.
\end{proof}

We will prove now,  that $\varphi^{-1}$ is a power of  the denominator of $x$,
multiplied by some constant from $k$.

\begin{lemma}\label{lem:relationOfuAndyinKtuAndKxyLinFactorIsDenominatPow}
  Let  $k(t,u)=k(x,y)$, $u^2=D_t$,  $y^2=D_x$  be  a hyperelliptic  function
  field
  % TODO: allow g=1 if lem:kxEqKtInHyperellFields is true in this case
  over  a  constant  field  $k$  of  characteristic  $\neq  2$,  where  both
  $D_t  \in  k[t]$ and  $D_x  \in  k[x]$  are monic  separable  polynomials.
  Let   $x=\frac{\alpha_0t+\alpha_1}{\alpha_2t+\alpha_3}$,   $\alpha_i   \in
  k$,    $\alpha_0\alpha_3-\alpha_1\alpha_2   \neq    0$   as    stated   in
  corollary~\ref{cor:relationOftAndxinKtuAndKxy}  and  $y=\varphi u$  as  in
  lemma~\ref{lem:relationOfuAndyinKtuAndKxyLinear}.  Then there  are $\gamma
  \in \unitGroup{k}$ and $m \in \N$ s.th.\
  \[
    \varphi^{-1} = \gamma (\alpha_2t+\alpha_3)^m.
  \]
\end{lemma}
\begin{proof}
  By  lemma~\ref{lem:relationOfuAndyinKtuAndKxyLinFactorIsPolyInv}  we  know
  $\varphi^{-1}  \in  k[t]\setminus \{0\}$.  Factoring  it  over $k$  yields
  $\varphi^{-1}   =   \gamma   \cdot  (\alpha_2t   +   \alpha_3)^m$,   where
  $\gamma   \in  k[t]\setminus   \{0\}$  s.th.\   $(\alpha_2t  +   \alpha_3)
  \ndivides  \gamma$  ($\gamma$  does  not   need  to  be  irreducible).  By
  lemma~\ref{lem:relationOfuAndyinKtuAndKxyFactorsOfDxWrtt} we have
  \[
    D_t = \varphi^{-2} (\alpha_2t+\alpha_3)^{-d_x} \prod_{i=1}^{d_x} p_i
        = \gamma^{2} (\alpha_2t+\alpha_3)^{m-d_x} \prod_{i=1}^{d_x} p_i.
  \]
  As $D_t$ is separable, we need to have $\gamma \in \unitGroup{k}$ which
  proves our claim.
\end{proof}

Computing the degree of $\varphi^{-1}$, we see that it is a scalar multiple of
the $(g+1)$-th power of the denominator of $x$.

\begin{lemma}\label{lem:relationOfuAndyinKtuAndKxyLinFactorFormula}
  Let  $k(t,u)=k(x,y)$,  $u^2=D_t$,  $y^2=D_x$  be  a  hyperelliptic  function
  field  over a  constant field  $k$ of  characteristic $\neq  2$, where  both
  $D_t  \in  k[t]$  and  $D_x  \in  k[x]$  are  monic  separable  polynomials.
  Let  $x =  \frac  {\alpha_0t +\alpha_1}  {\alpha_2t +\alpha_3}$,  $y=\varphi
  u$   as   stated   in   corollary~\ref{cor:relationOftAndxinKtuAndKxy}   and
  lemma~\ref{lem:relationOfuAndyinKtuAndKxyLinear}. Then we have
  \begin{enumerate}
  \item
    If $x \in k[t]$, then $\varphi \in \unitGroup{k}$.
  \item
    If $x  \notin k[t]$, we  assume w.l.o.g.\ $\alpha_2=1$. Then  there exists
    some  $\gamma \in  \unitGroup{k}$  s.th.\ $\varphi^{-1}  = \gamma  (t +
    \alpha_3)^{g+1}$.
  \end{enumerate}
\end{lemma}
\begin{proof}
  By   lemma~\ref{lem:relationOfuAndyinKtuAndKxyFactorsOfDxWrtt},  there   are
  $p_i    \in   \AlgebraicClosure{k}[t]$,    s.th.\   $D_t    =   \varphi^{-2}
  (\alpha_2t  +   \alpha_3)^{-d_x}  \prod_{i=1}^{d_x}  p_i$  and   $d_x-1  \le
  \deg_t\left(\prod_{i=1}^{d_x}  p_i \right)  \le  d_x$. Let  us consider  the
  given cases, separately.
  \begin{enumerate}
  \item
    Let     us     assume     $x\in     k[t]$,     first.     We     already
    know       $\varphi^{-1}      \in       k[t]      \setminus       \{0\}$
    (cf.\  lemma~\ref{lem:relationOfuAndyinKtuAndKxyLinFactorIsPolyInv}) and
    $D_t=\varphi^{-2}\alpha_3^{-d_x}\prod_{i=1}^{d_x} p_i$. If $\varphi^{-1}
    \notin k$, then  $\varphi^{-2}$ were a non trivial  square polynomial in
    $t$ dividing  $D_t$. This  contradicts the  separability of  $D_t$. Thus
    $\varphi^{-1} \in k$, which immediately implies $\varphi \in \unitGroup{k}$.
  \item
    We  proceed with  the case  $x  \notin  k[t]$,  i.e.\  $\alpha_2 \neq  0$.
    By  reducing  the  fraction  $x= \frac  {\alpha_0t  +\alpha_1}  {\alpha_2t
    +\alpha_3}$, we  can assume  $\alpha_2=1$ without  loss of  generality. As
    $\varphi^{-1} \in k[t]$, we get
    \begin{align*}
       \deg_t(D_t)
         =& 2\deg_t(\varphi^{-1})
         - d_x\deg_t(t+\alpha_3)
         + \deg_t\left(\prod_{i=1}^{d_x} p_i \right)
         \\
         =& 2\deg_t(\varphi^{-1})
         - d_x
         + \deg_t\left(\prod_{i=1}^{d_x} p_i \right),
    \end{align*}
    which implies
    \[
       2\deg_t(\varphi^{-1})
        =
         \deg_t(D_t)
       + d_x
       - \deg_t\left(\prod_{i=1}^{d_x} p_i \right).
    \]
    Thus,        the        inequality        $d_x-1        \le
    \deg_t\left(\prod_{i=1}^{d_x} p_i \right) \le d_x$ yields
    \begin{align*}
        \deg_t(D_t)
        =&
         \deg_t(D_t)
       + d_x
        -d_x
         \\
      \le &
         \deg_t(D_t)+d_x-\deg_t\left(\prod_{i=1}^{d_x} p_i \right)
         \\
      = &
       2 \deg_t(\varphi^{-1})
       \\
      \le &
         \deg_t(D_t)
       + d_x
         -d_x +1
       \\
       =&
       \deg_t(D_t) + 1.
    \end{align*}
    As      $\deg_t      (D_t)      \in     \{      2g+1,      2g+2      \}$
    we     conclude     $\deg_t      (\varphi^{-1})     =     g+1$.     From
    lemma~\ref{lem:relationOfuAndyinKtuAndKxyLinFactorIsDenominatPow}     we
    know  that there  is  some $\gamma  \in \unitGroup{k}$  and some  $m
    \in  \N$ s.th.\  $\varphi^{-1} =  \gamma (t  + \alpha_3)^m$.  As
    $\deg_t(\varphi^{-1})=g+1$,  this  implies  our claim.
  \end{enumerate}
\end{proof}

\subsection{Putting Both Relations Together}
  The following theorem completely characterizes  the relation between any two
  bases of a hyperelliptic function field of characteristic $\neq 2$. This can
  be used to  check whether a  given function field has a  specific kind of
  defining equation. It is the key ingredient of algorithm~\ref{alg:compAutWrtAlgClosedK},
  which computes
  automorphism groups.

  Using the  facts proved above,  it remains to  compute the scalar  factor of
  $\varphi$ in order to know the relation between two bases:

  \begin{theorem}\label{thm:relationOfGeneratorsOfHyperellipticFields}
    Let    $k(t,u)=k(x,y)$,   $u^2=D_t$,    $y^2=D_x$   be    a
    hyperelliptic function  field
    % TODO: allow g=1 if lem:kxEqKtInHyperellFields is true in this case
    over a constant
    field $k$  of characteristic  $\neq 2$,  where both  $D_t \in
    k[t]$ and $D_x \in k[x]$ are monic separable polynomials. Let
    $d_x:=\deg_x(D_x)$.
    \begin{enumerate}
    \item
      If  $x \in  k[t]$,  then there  are $\alpha_0,  \alpha_1  \in k$  s.th.\
      $x=\alpha_0t+\alpha_1$  and  $\alpha_0  \neq  0$.  Furthermore  we  have
      $y=\varphi u$ with $\varphi \in \unitGroup{k}$,
      \[
         \varphi^2 = \alpha_0^{d_x}.
      \]
    \item
      If  $x  \notin  k[t]$,  then there  are  $\alpha_0,  \alpha_1,  \alpha_3
      \in  k$,  s.th.\  $x=\frac{\alpha_0t  +  \alpha_1}{t  +  \alpha_3}$,
      \mbox{$\alpha_0\alpha_3-\alpha_1   \neq   0$}.   Furthermore   we   have
      $y=\varphi u$, where
      \[
        \varphi = \frac{\beta} {(t + \alpha_3)^{g+1}},
      \]
      with $\beta \in \unitGroup{k}$. For $\beta$ we have the formula
      \[
        \beta^2 = \begin{cases}
                     D_x (\alpha_0 )
                      &\text{,\ if\ } D_x (\alpha_0) \neq 0
                    \\
                     (\alpha_1 - \alpha_0 \alpha_3)
                      \tilde D_x (\alpha_0)
                      &\text{,\ if\ } D_x (\alpha_0) = 0,
                  \end{cases}
      \]
      where  $\tilde   D_x(x)  :=  \frac  {D_x(x)}   {x  -
      \alpha_0}$.
    \end{enumerate}
  \end{theorem}
  \begin{proof}
    Corollary~\ref{cor:relationOftAndxinKtuAndKxy}    gives   the    existence
    of    $\alpha_0, \dots, \alpha_3   \in   k$ s.th.\
    $x= \frac {\alpha_0t +\alpha_1} {\alpha_2t +\alpha_3}$ and $\alpha_0\alpha_3 -
    \alpha_1\alpha_2 \neq 0$.
    Lemma~\ref{lem:relationOfuAndyinKtuAndKxyLinear}   yields  some   $\varphi
    \in   k(t)    \setminus   \{0\}$    s.th.\   $y    =   \varphi    u$.   By
    lemma~\ref{lem:relationOfuAndyinKtuAndKxyLinFactorFormula},     we    know
    $\varphi  \in  \unitGroup{k}$   if  $x  \in  k[t]$   and  $\varphi^{-1}  =
    \gamma  (\alpha_2 t  +  \alpha_3)^{g+1}$ with  $\gamma \in  \unitGroup{k}$
    otherwise.
    Lemma~\ref{lem:relationOfuAndyinKtuAndKxyFactorsOfDxWrtt} implies
    \begin{equation}\label{eqn:DtFactorizationFromDxtSubstitution}
      D_t = \varphi^{-2} (\alpha_2t+\alpha_3)^{-d_x} \prod_{i=1}^{d_x} p_i,
    \end{equation}
    where $p_i =  (\alpha_0 - \alpha_2\eta_i)t + \alpha_1  - \alpha_3\eta_i$
    and the $\eta_i \in \AlgebraicClosure{k}$ are the zeroes of $D_x$.

    Let us consider the different cases, now:
    \begin{enumerate}
    \item
      If    $x   \in    k[t]$,    we   have    $\alpha_2=0$.   Reducing    the
      fraction    $\frac   {\alpha_0t    +\alpha_1}   {\alpha_3}$,    we   may
      w.l.o.g.\ assume    $\alpha_3=1$.    Thus
      equation~(\ref{eqn:DtFactorizationFromDxtSubstitution}) becomes
      \[
        D_t = \varphi^{-2} 
	   \prod_{i=1}^{d_x} (\alpha_0 t + \alpha_1 - \eta_i).
      \]
      As     $\alpha_0     \neq     0$     (which     we     conclude     from
      $\alpha_0\alpha_3-\alpha_1\alpha_2  =  \alpha_0   \neq  0$)  and
      $\varphi \in k$, the leading coefficient of $D_t$ is
      \[
        1 = \leadingCoeff[t]{D_t} = \varphi^{-2} \alpha_0^{d_x},
      \]
      because  $D_t$   is  monic  by  assumption.   This  implies
      $\varphi^2   =    {\alpha_0}^{d_x}$.
    \item
      If   $x   \notin   k[t]$,   we  have   $\alpha_2   \neq   0$.   Reducing
      the  fraction  $\frac   {\alpha_0t  +\alpha_1}  {\alpha_2t  +\alpha_3}$,
      we   may    assume   $\alpha_2=1$. We already know  $\varphi^{-1}    =   \gamma
      (t    +     \alpha_3)^{g+1}$. Setting    $\beta:=\gamma^{-1}$,
      it remains to compute $\beta^2$.
      From equation~(\ref{eqn:DtFactorizationFromDxtSubstitution}), we get
      \begin{equation*}
          D_t =
            \beta^{-2} (t+\alpha_3)^{2g+2-d_x} \prod_{i=1}^{d_x} p_i.
      \end{equation*}
      As before, we compute the leading coefficients:
      \begin{equation*}
         1 = \leadingCoeff[t]{D_t} 
	 = \leadingCoeffLarge[t] {\beta^{-2} (t+\alpha_3)^{2g+2-d_x} \prod_{i=1}^{d_x} p_i} 
	 = \beta^{-2} \leadingCoeffLarge[t] {\prod_{i=1}^{d_x} p_i}.
      \end{equation*}
      We obtain
      \begin{equation}\label{eqn:numeratorOfVarphiIsLeadingCoeffOfDxtSubstitution}
        \beta^2 = \leadingCoeffLarge[t] {\prod_{i=1}^{d_x} p_i}.
      \end{equation}
      From    Lemma~\ref{lem:relationOfuAndyinKtuAndKxyFactorsOfDxWrtt},    we
      know  $d_x   -1  \le  \deg_t(\prod_{i=1}^{d_x}  p_i)   \le  d_x$.  Thus,
      there  are   two  cases:  $\deg_t(\prod_{i=1}^{d_x}  p_i)   =  d_x$  and
      $\deg_t(\prod_{i=1}^{d_x}   p_i)  =   d_x-1$.   In   the  latter   case,
      there  is   some  index  $j$   s.th.\  $p_j  =  (\alpha_0   -\eta_j)t  +
      \alpha_1-\alpha_3\eta_j  \in k$,  i.e.\ $\alpha_0-  \eta_j=0$. Hence,
      ${\alpha_0}  = \eta_j$, which implies  $D_x({\alpha_0}) =
      0$.  In  the  former  case,  there  is no  such  index,  i.e.\  we  have
      $D_x({\alpha_0}) \neq 0$.

      \begin{enumerate}
      \item
        If  $D_x({\alpha_0}) \neq  0$, we  have
        $\alpha_0  - \eta_i  \neq 0$  for all  $i$. Thus  we get
        \begin{align*}
          \beta^2 
	  =& \leadingCoeffLarge[t] {\prod_{i=1}^{d_x} p_i}
          = \leadingCoeffLarge[t] {\prod_{i=1}^{d_x} (\alpha_0 - \eta_i)t 
	        + \alpha_1 -\alpha_3\eta_i}
	  \\
          =& \prod_{i=1}^{d_x} (\alpha_0 - \eta_i)
          = D_x({\alpha_0})
        \end{align*}
        as claimed.
      \item
        If $D_x({\alpha_0}) = 0$, there is 
        exactly one index $j$ s.th.\ $\eta_j = {\alpha_0}$.
        W.l.o.g.\ we assume $j=d_x$. Then we have
        $p_{d_x} = \alpha_1 - \alpha_0\alpha_3$.
        Thus equation~(\ref{eqn:numeratorOfVarphiIsLeadingCoeffOfDxtSubstitution})
        implies
        \begin{align*}
          \beta^2 
	  =& \leadingCoeffLarge[t] {\prod_{i=1}^{d_x} p_i}
          = \leadingCoeffLarge[t] {\prod_{i=1}^{d_x} (\alpha_0 - \eta_i)t 
	        + \alpha_1 -\alpha_3\eta_i}
	  \\
          =& \leadingCoeffLarge[t] {(\alpha_1 - \alpha_0\alpha_3)
	      \prod_{i=1}^{d_x-1} (\alpha_0 - \eta_i)t + \alpha_1 -\alpha_3\eta_i}
	  \\
          =& (\alpha_1 - \alpha_0\alpha_3) \prod_{i=1}^{d_x-1} (\alpha_0 - \eta_i)
          \\
          =& (\alpha_1 - {\alpha_0\alpha_3})
                \tilde D_x ({\alpha_0}),
        \end{align*}
        as   $\tilde  D_x   (x)  =   \frac  {D_x   (x)}  {x-
        {\alpha_0}} = \frac {D_x (x)} {x-\eta_{d_x}} =
        \prod_{i=1}^{d_x-1} (x  - \eta_i)$. 
      \end{enumerate}
    \end{enumerate}
  \end{proof}

\begin{corollary}\label{cor:relationOfGeneratorsOfHyperellipticFieldsIffTwoBases}
  Let  $k(x,y)$,   $y^2=D_x$  be  a   hyperelliptic  function  field   over  a
  constant   field  $k$   of   characteristic  $\neq   2$,   where  $D_x   \in
  k[x]$   is  a   monic  separable   polynomial.   Let  $D_t   \in  k[T]$   be
  another  monic  separable   polynomial.  There  exists  a   basis  $t,u  \in
  k(x,y)$   s.th.\  $k(x,y)=k(t,u)$,   $u^2=D_t(t)$  iff   there  exist   $t,u
  \in  k(x,y)$   for  which  $u^2=D_t(t)$   and  the  relations  $x   =  \frac
  {\alpha_0t  +\alpha_1} {\alpha_2t  +\alpha_3}$,  $y =  \varphi  u$ given  in
  theorem~\ref{thm:relationOfGeneratorsOfHyperellipticFields} hold.
\end{corollary}
\begin{proof}
  It remains to  show that the existence  of $t, u$, $u^2=D_t(t)$  s.th.\ $x =
  \frac  {\alpha_0t +\alpha_1}  {\alpha_2t  +\alpha_3}$, $y  =  \varphi u$  as
  given in theorem~\ref{thm:relationOfGeneratorsOfHyperellipticFields} implies
  $k(x,y)=k(t,u)$. It  is obvious, that  $k(x) \subseteq k(t)$ and $k(t)(u)=k(t)(y)$.
  Solving $x = \frac  {\alpha_0t +\alpha_1}  {\alpha_2t  +\alpha_3}$
  for  $t$, 
  we see  $k(t)  \subseteq k(x)$.  Thus
  $k(t)=k(x)$,  i.e.\  $k(t,u) =  k(t)(u) =
  k(t)(y) = k(x)(y)=k(x,y)$.
\end{proof}

%\begin{remark}
%  Without   additional    assumptions,   there    exist   no    formulas   for
%  the    relationship     between    two    bases    of     a    hyperelliptic
%  function   field   that   are   more    specific   than   those   given   in
%  theorem~\ref{thm:relationOfGeneratorsOfHyperellipticFields}. 
%  
%  Let  us   prove  this   statement:  First   of  all,   for  each   choice  of
%  $\alpha_0,  \alpha_1, \alpha_2,  \alpha_3  \in  k$ s.th.\  $\alpha_0\alpha_3
%  -\alpha_1\alpha_2  \neq  0$,   it  is  well  known  that   $k(t)  =  k\left(
%  \frac   {\alpha_0t    +\alpha_1}   {\alpha_2t    +\alpha_3}\right)$,   since
%  $\AutomorphismGroup{k(t)}{k} = \PGL{2}{k}$. Furthermore,  we can reduce this
%  fraction to either $\alpha_0't+\alpha_1'$  or $\frac {\alpha_0't +\alpha_1'}
%  {t +\alpha_3}$,  yielding the  relation between  $t$ and  $x$ given  in our
%  theorem.  The relation  $y =  \varphi u$  can also  not be  ``improved'': We
%  already know  $\varphi$ up to its  sign, and its negative  value yields $-y=
%  (-\varphi)u$. As $(-y)^2 = y^2 = D_x$,  this gives a basis $x,-y$ which
%  cannot be distinguished from $x,y$ by the assumptions of our theorem.
%
%  Of    course,   our    relations   can    be   simplified,    if   additional
%  assumptions   are    made.   One    such   theorem    can   be    found   in
%  \cite[Proposition~1.2]{lockhart:1994}. Unfortunately,  Lockhart's assumptions
%  do not hold in our situation.
%\end{remark}
%

As  said  before,  theorem~\ref{thm:relationOfGeneratorsOfHyperellipticFields}
can  be  applied   to  check  if  a  hyperelliptic   function  field  $k(x,y)$
has  a   basis  $t,u$  satisfying   a  given  equation   $u^2=D_t$.  According
to   corollary~\ref{cor:relationOfGeneratorsOfHyperellipticFieldsIffTwoBases},
we    can    decide   this    question    by    checking,   if    there    are
$\alpha_i    \in    k$,    and    $\varphi   \in    k(t)$    as    given    in
theorem~\ref{thm:relationOfGeneratorsOfHyperellipticFields}   s.th.\  $u^2   =
D_t$, which  is equivalent  to $D_t  = u^2 =  \varphi^{-2} y^2  = \varphi^{-2}
D_x$, here. This can be done using the following algorithm:
\begin{algorithm}\label{alg:defEqWrtK}
  Let  $k(x,y)$,  $y^2=D_x$,  $D_x\in  k[x]$  monic  and  separable,  be  some
  hyperelliptic function field  of genus $g$ with  $\characteristic{k} \neq 2$
  and let $D_t \in  k[t]$ be some monic, separable polynomial  of $\deg_t(D_t) \in
  \{ 2g+1, 2g+2 \}$. Let $d_x:=\deg_x(D_x)$.
  \begin{enumerate}
  \item
    We compute $\varphi^2 \in k(t)$ symbolically from the $\alpha_i$ according
    to  theorem~\ref{thm:relationOfGeneratorsOfHyperellipticFields}. Since  we
    do not know the  $\alpha_i$ in advance, we cannot tell  which of the cases
    of our  theorem applies. Thus  we have to  compute $\varphi^2$ and  do the
    following steps in each of these cases:
    \begin{itemize}
    \item
      If     $x\in     k[t]$,     we     have to use
      $x=\alpha_0t     +\alpha_1$,
      $\varphi^2=\alpha_0^{d_x}$.
    \item
      If $x\notin k[t]$, i.e.\  $x= \frac {\alpha_0t +\alpha_1} {t+\alpha_3}$,
      we  consider both  $D_x(\alpha_0)\neq 0$  and $D_x(\alpha_0)=0$.  In the
      former case we have $\varphi^2 = D_x(\alpha_0) (t+\alpha_3)^{-2g+2}$. If
      $D_x(\alpha_0) =  0$, we know that  $x-\alpha_0$ is a divisor  of $D_x$.
      Thus  we  can  find  all possible  $\alpha_0$  explicitly  by  factoring
      $D_x$ over  $k$. For  each such  $\alpha_0$, we  compute $\tilde  D_x :=
      \frac {D_x(x)}{x-\alpha_0}$ obtaining $\varphi^2  = (\alpha_1 - \alpha_0
      \alpha_3) \tilde D_x (\alpha_0) (t+\alpha_3)^{-2g-2}$.
    \end{itemize}
  \item  
    After multiplying by  the denominators, our condition  $D_t = \varphi^{-2}
    D_x$ becomes  an equation  of polynomials  in $t$  and the  $\alpha_i$. We
    compare coefficients of $t$. The  resulting system of polynomial equations
    for the $\alpha_i$ is denoted by $(*)$.
  \item
    Let   the  ideal   $I$  be   generated   by  $(*)$   and  the   polynomial
    $1-(\alpha_0\alpha_3  - \alpha_1\alpha_2)T$,  where  $T$  is new  variable
    symbol  and  the  $\alpha_i$  satisfy $x=  \frac  {\alpha_0t  +  \alpha_1}
    {\alpha_2t +\alpha_3}$  according to  the case  we are  considering. Using
    Gr\"obner  basis methods,  we check  $I$ for  solvability and  construct a
    solution, if it exists.
  \end{enumerate}
  Thus we can construct a basis $k(x,y)=k(t,u)$, $u^2=D_t$ 
  iff there are $\alpha_i, T$ in the variety of $I$ over $k$
  for any of the cases mentioned in step (1).
\end{algorithm}

Let us illustrate this algorithm with an example: 
\begin{example}\label{ex:findingABasisTransformation}
  Let $k=\GF{11}$, $F  = k(x,y)$, $y^2 =  D_x := x^5+x^4+4x^3+5x^2+10x+7$.
  We would like to  know, if there is a basis $F = k(t,u)$  s.th.\ $u^2 = D_t
  := t^5+7t^3+9t^2+9t+6$.

  We    start    with    the    easiest    case    $x\in    k[t]$.    From
  theorem~\ref{thm:relationOfGeneratorsOfHyperellipticFields},    we   get
  $x=\alpha_0t+ \alpha_1$, $y=\varphi u$ and $\varphi^2 = \alpha_0^5$.

  Substituting, we get
  \begin{align*}
     D_x= D_x(\alpha_0t +\alpha_1)=&
     \alpha_0^5 t^5\\
     &+(5\alpha_0^4\alpha_1+\alpha_0^4)t^4 \\
     &+(10\alpha_0^3\alpha_1^2+4\alpha_0^3\alpha_1+4\alpha_0^3)t^3 \\
     &+(10\alpha_0^2\alpha_1^3+6\alpha_0^2\alpha_1^2+\alpha_0^2\alpha_1+5\alpha_0^2)t^2 \\
     &+(5\alpha_0\alpha_1^4+4\alpha_0\alpha_1^3+\alpha_0\alpha_1^2+10\alpha_0\alpha_1
       +10\alpha_0)t \\
       &+\alpha_1^5+\alpha_1^4+4\alpha_1^3+5\alpha_1^2+10\alpha_1+7
  \end{align*}
  where $\alpha_0,  \alpha_1$ are to  be found. Comparing  coefficients in
  $D_x= \varphi^2 D_t = \alpha_0^5 D_t$ yields the equations $(*)$:
  \begin{align*}
     5\alpha_0^4\alpha_1+\alpha_0^4&=0, \\
     10\alpha_0^3\alpha_1^2+4\alpha_0^3\alpha_1+4\alpha_0^3&= 7\alpha_0^5, \\
     10\alpha_0^2\alpha_1^3+6\alpha_0^2\alpha_1^2+\alpha_0^2\alpha_1+5\alpha_0^2
        &=9\alpha_0^5, \\
     5\alpha_0\alpha_1^4+4\alpha_0\alpha_1^3+\alpha_0\alpha_1^2+10\alpha_0\alpha_1+10
        \alpha_0&=9\alpha_0^5, \\
     \alpha_1^5+\alpha_1^4+4\alpha_1^3+5\alpha_1^2+10\alpha_1+7 &= 6\alpha_0^5
  \end{align*}
  Augmenting  $(*)$  by   $1-\alpha_0T$,  we  get  the   ideal  $I$.  Singular
  (\cite{SINGULAR}) computes the following  Gr\"obner basis of $I$ w.r.t.\
  the lexicographical ordering:
  \begin{align*}
    T-4 =&0 \\
    \alpha_0-3 =&0 \\
    \alpha_1-2\alpha_0^5T^2-3\alpha_0^2T^2 =&0
  \end{align*}
  This  implies  $T=4$,  $\alpha_0=3$.  Substituting  these  values  into  the
  remaining equation, we obtain $\alpha_1=2$. Thus, setting $t:=4x-3$, $u:=y$,
  i.e.\ $x=3t+2$, $\varphi=3^5=1$, we get a basis $F=k(t,u)$, with $u^2=D_t$.
\end{example}

In   order    to   compute   the   automorphism    group   $\AutomorphismGroup
{\AlgebraicClosure{k}  (x,y)}  {\AlgebraicClosure{k}}$ over  an  algebraically
closed constant field,  it suffices to check  $\AlgebraicClosure{k} (x,y)$ for
normal forms, as we will see in section~\ref{sec:compAut}. This simplifies the
Gr\"obner basis  step of  algorithm~\ref{alg:defEqWrtK}, giving  the following
modified algorithm:

\begin{algorithm}\label{alg:defEqWrtAlgClosureOfK}
  Let  $k(x,y)$,  $y^2=D_x$,  $D_x\in  k[x]$  monic  and  separable,  be  some
  hyperelliptic function field  of genus $g$ with  $\characteristic{k} \neq 2$
  and $D_t \in  k[t]$ be some monic, separable polynomial  of $\deg_t(D_t) \in
  \{ 2g+1, 2g+2 \}$. Let $d_x:=\deg_x(D_x)$.

  Whether  there     exists     a     basis     $\AlgebraicClosure{k}(x,y)     =
  \AlgebraicClosure{k}(t,u)$  with  $u^2=D_t$,  can be  checked  analogous  to
  algorithm~\ref{alg:defEqWrtK}. We only note the differences:
  \begin{enumerate}
  \item
    In  order   to  compute   $\varphi^2$  in  the   case  $x   \notin  k[t]$,
    $D_x(\alpha_0)=0$, we have to consider all zeroes $\alpha_0$ of $D_x$ over
    $\AlgebraicClosure{k}$, i.e.\ we  have to factor $D_x$  over its splitting
    field.
  \item[(3)]
    Instead of constructing an element of the  variety of $I$, we only need to
    check if  it's empty. To do  so, we compute  a Gr\"obner basis $B$  of $I$
    (e.g.\ w.r.t.\ the degree  reverse lexicographical ordering). There exists
    a solution  $\alpha_i, T  \in \AlgebraicClosure{k}$,  iff $I  \neq \langle
    1\rangle$, i.e.\ iff $B \neq \{1\}$.
  \end{enumerate}
  As  in algorithm~\ref{alg:defEqWrtK},  we  infer the  existence  of a  basis
  $\AlgebraicClosure{k}(x,y)   =  \AlgebraicClosure{k}(t,u)$,   $u^2=D_t$  iff
  $B\neq \{1\}$.
\end{algorithm}

\begin{remark}
  An       essential      feature       of      algorithms~\ref{alg:defEqWrtK}
  and~\ref{alg:defEqWrtAlgClosureOfK} is, that $D_t$ does not need to be known
  completely. It  may contain  some parameters  for which  we can  also solve.
  Therefore,  we can  use our  algorithms to  check, whether  a given  hyperelliptic
  function field $\AlgebraicClosure{k}(x,y)$ has some of Brandt's normal forms
  (cf.\ theorem~\ref{thm:brandtsNormalForms}). We will see  how to do this, in
  the following section.
\end{remark}

\section{Computing the Automorphism Group}\label{sec:compAut} 
\subsection{Algebraically Closed Constant Fields}

Algorithm~\ref{alg:defEqWrtAlgClosureOfK}  can  be   applied  to  compute  the
automorphism group  of a  hyperelliptic function  field over  an algebraically
closed constant  field: 

\begin{algorithm}\label{alg:compAutWrtAlgClosedK}
  Let  $k(x,y)$,  $y^2  =  D_x$,  $D_x  \in k[x]$  monic  and  separable,  be  a
  hyperelliptic function field of genus $g$ and $\characteristic{k}\neq 2$.
  We denote $F:= \AlgebraicClosure{k} (x,y)$.
  \begin{enumerate}
  \item
    For   each  possible   type  $\typeOfField   {G}  {\AlgebraicClosure{k}}$,
    we    look   up    the    corresponding   normal    form   $u^2=D_t$    in
    theorem~\ref{thm:brandtsNormalForms}.
  \item  
    For each normal form found in step (1), we check,  
    for which  parameter sets  $D_t$ 
    has  degree $2g+1$  or $2g+2$. 
    This yields the set $N$ of all polynomials $D_t$, s.th.\ $u^2=D_t$
    is  a  normal  form for  a  field  of  genus  $g$ and  type  $\typeOfField
    {G}  {\AlgebraicClosure{k}}$. The integer  parameters in  each $D_t
    \in  N$ are  fixed,  while the  $D_t$ may  still  contain parameters  from
    $\AlgebraicClosure{k}$.
  \item
    For each $G$, $N$ and each $D_t \in N$, 
    we check if $F$ has a basis $F = \AlgebraicClosure{k} (t,u)$
    satisfying $u^2 = D_t$ as well as the additional conditions from
    theorem~\ref{thm:brandtsNormalForms}.
    To do so, we use a slight 
    modification of algorithm~\ref{alg:defEqWrtAlgClosureOfK}:
    
    Let $C_{0}$ and $C_{1}$ be the sets of polynomials that
    according to theorem~\ref{thm:brandtsNormalForms} have to be $=0$
    and $\neq 0$, respectively. 
    Let
    \[
      c := (\alpha_0\alpha_3 -\alpha_1\alpha_2)
         \prod_{f \in C_{1}} f.
    \]
    We define the ideal $I$ to be generated by $(*)$, $C_{0}$ and $1-c\cdot T$
    rather than  just by $(*)$ and  $1-(\alpha_0\alpha_3 -\alpha_1\alpha_2)T$.
    Note that the  polynomial ring $R \supseteq I$ may  contain more variables
    than just the $\alpha_i$ and $T$, now.

    We apply the rest of algorithm~\ref{alg:defEqWrtAlgClosureOfK} without any
    changes.
    
    The variety of $I$ is non-empty iff $k(x,y)$ is of type $\typeOfField {G}
    {\AlgebraicClosure{k}}$.
  \item
    Let $G$  be the largest group  $G$ s.th.\ $k(x,y)$ is  of type $\typeOfField
    {G} {\AlgebraicClosure{k}}$. Then
    \[
      \AutomorphismGroup {\AlgebraicClosure{k} (x,y)} {\AlgebraicClosure{k}}  
	/ \CyclicGroup{2} \cong G,
    \]
    and  the generators  of  $\typeOfFieldAutSubGroup  {G} =  \AutomorphismGroup
    {\AlgebraicClosure{k}   (x,y)}   {\AlgebraicClosure{k}}$    are   given   in
    theorem~\ref{thm:brandtsNormalForms}.
  \end{enumerate}
\end{algorithm}

Thus,  we  are  able to  compute  the  structure  as  well as  the  generators
of  the automorphism  group  $\AutomorphismGroup {\AlgebraicClosure{k}  (x,y)}
{\AlgebraicClosure{k}}$ for each hyperelliptic function field $k(x,y)$.

\begin{example}\label{ex:computingAutOverAlgClosure}
  Let  $F:=\AlgebraicClosure{\GF{7}}(x,y)$   with  $y^2=  x^5+x^3+x$   as  in
  example~\ref{ex:brandtsNormalFormPartiallyVisible}.   The  above   algorithm
  yields that $F$ is of the types
  $\typeOfField{\CyclicGroup{2}}{\AlgebraicClosure{\GF{7}}}$,
  $\typeOfField{\CyclicGroup{3}}{\AlgebraicClosure{\GF{7}}}$,
  $\typeOfField{\CyclicGroup{6}}{\AlgebraicClosure{\GF{7}}}$,
  $\typeOfField{\DihedralGroup{2}}{\AlgebraicClosure{\GF{7}}}$,
  $\typeOfField{\DihedralGroup{3}}{\AlgebraicClosure{\GF{7}}}$ and
  $\typeOfField{\DihedralGroup{6}}{\AlgebraicClosure{\GF{7}}}$.
  
  To see, how the algorithm works, we  consider parts of the proof that $F$ is
  of type $\typeOfField{\DihedralGroup{3}}{\AlgebraicClosure{\GF{7}}}$:
  \begin{enumerate}
  \item 
    The        normal       form        for        fields       of        type
    $\typeOfField{\DihedralGroup{3}}{\AlgebraicClosure{\GF{7}}}$  is given  by
    \[
      y^2=t^{\nu_0}  (t^3-1)^{\nu_1}   (t^3+1)^{\nu_2}  \prod_{i=1}^s   (t^6  -
      a_jt^3+1),
    \]
    where  $\nu_i \in \{0,1\}$,  $\nu_1=\nu_2$, $s\in \N$  and the
    $a_j \in \AlgebraicClosure{\GF{7}}\setminus \{\pm 2\}$ are pairwise distinct.
  \item
    As $g=2$, we need to have $\deg_t(D_t) \in \{ 5, 6\}$, from which we get
    \[
      N = \{ (t^3 - 1)(t^3 + 1) , t^6 - a_1t^3 + 1 \}
    \]
  \item
    Algorithm~\ref{alg:defEqWrtAlgClosureOfK}  finds    out    
    that    $F$    possesses    a    basis
    $F=\AlgebraicClosure{\GF{7}}(t,u)$,  $u^2=(t^3-1)(t^3+1)$.  Thus,  $F$  is
    of    type   $\typeOfField{\DihedralGroup{3}}{\AlgebraicClosure{\GF{7}}}$.
    The   corresponding   system   $(*)$   of   equations   and   inequalities
    is   not   given,  as   it   looks   quite   ugly   and  does   not   help
    in   understanding   this  step   of   the   algorithm.  It   is   similar
    to   the   one  given   in   example~\ref{ex:findingABasisTransformation}.
    Simplifying $(t^3-1)(t^3+1) = t^6-1$, theorem~\ref{thm:brandtsNormalForms}
    immediately       implies       that       $F$      is       of       type
    $\typeOfField{\DihedralGroup{6}}{\AlgebraicClosure{\GF{7}}}$.

    Furthermore, the second  element of $N$ can  also be used to  find a basis
    of  $F$:  Setting  $a_1:=0$  and $\alpha_1:=i$,  where  $i^2=-1$,  implies
    $\alpha_0=1$  and  $\alpha_3=-i$. Thus  $F=\AlgebraicClosure{\GF{7}}(v,w)$
    with $w^2 = v^6+1$.
  \item
    From    the    list    above,    we    know    that    $\DihedralGroup{6}$
    is    the    largest    group    $G$     s.th.\    $F$    is    of    type
    $\typeOfField {G}  {\AlgebraicClosure{\GF{7}}}$. Thus, $\AutomorphismGroup
    {\AlgebraicClosure{\GF{7}}(x,y)}    {\AlgebraicClosure{\GF{7}}}$   is    a
    central  extension   of  $\DihedralGroup{6}$  by   the  $\CyclicGroup{2}$,
    generated   by    the   hyperelliptic   involution.   From    the   normal
    form   $u^2   =   (t^6+1)$,    we   know   $\nu_0=\nu_1=s=0$,   $\nu_2=1$.
    According    to    theorem~\ref{thm:brandtsNormalForms},    a    set    of
    generators    of   $\AutomorphismGroup    {\AlgebraicClosure{\GF{7}}(x,y)}
    {\AlgebraicClosure{\GF{7}}}$  is  given  by $\{\varphi,  \psi,  \sigma\}$,
    where $\varphi: t  \mapsto t$, $u \mapsto -u$, $\psi:  t \mapsto \zeta t$,
    $u  \mapsto u$  and $\sigma:  t \mapsto  \frac 1t$,  $u \mapsto  \frac {u}
    {t^3}$ with a primitive $6$-th root $\zeta$ of unity.

    Looking   at    these   generators,   we    conclude   $\AutomorphismGroup
    {\AlgebraicClosure{\GF{7}}(x,y)}     {\AlgebraicClosure{\GF{7}}}     \cong
    \DihedralGroup{6} \times \CyclicGroup{2}$.
  \end{enumerate}
\end{example}

\subsection{Arbitrary Constant Fields}\label{subsec:compAut:arbitraryConstants}

Using algorithm~\ref{alg:compAutWrtAlgClosedK}, it is also possible to compute
$\AutomorphismGroup {k(x,y)}{k}$ for a  hyperelliptic function field $k(x,y)$,
where $k$ needs  not to be algebraically closed. A  similar application is the
computation  of  the smallest  algebraic  extension  $k' \supseteq  k$  s.th.\
$\AutomorphismGroup {k'(x,y)} {k'}  = \AutomorphismGroup {\AlgebraicClosure{k}
(x,y)} {\AlgebraicClosure{k}}$.

Let $k$ be  any field of characteristic $p>2$ and  $k(x,y)$ be a hyperelliptic
function  field. We  use  algorithm~\ref{alg:compAutWrtAlgClosedK} to  compute
the   types  of   $\AlgebraicClosure{k}(x,y)$.   Let  $k(x,y)$   be  of   type
$\typeOfField {G} {\AlgebraicClosure  {k}}$ and let $\AlgebraicClosure{k}(x,y)
=  \AlgebraicClosure{k}(t,u)$, $u^2=D_t$  be  the  corresponding normal  form.
Solving  the ideal  $I$ for  the $\alpha_i$  and the  parameters of  $D_t$, we
obtain explicit formulas for  the generators of $\typeOfFieldAutSubGroup {G}$.
Using these, it  is easy to find  the smallest field $k'  \supseteq k$, s.th.\
all automorphisms  from $\typeOfFieldAutSubGroup {G}$ define  automorphisms of
$k'(x,y)$.  Then, $k'  \supseteq  k$  is the  smallest  field extension  s.th.
$k(x,y)$ is of type $\typeOfField {G} {k'}$.

This method is used to solve the two problems given above: In order to compute
$\AutomorphismGroup  {k(x,y)} {k}$,  we  construct $k'$  for  each $G$  s.th.\
$k(x,y)$ is of type $\typeOfField {G} {\AlgebraicClosure{k}}$. The largest $G$
with $k'=k$ yields $\typeOfFieldAutSubGroup  {G} = \AutomorphismGroup {k(x,y)}
{k}$.

To find  the smallest  $k' \supseteq  k$ s.th.\  $\AutomorphismGroup {k'(x,y)}
{k'} =  \AutomorphismGroup {\AlgebraicClosure  {k} (x,  y)} {\AlgebraicClosure
{k}}$,  we   compute  $\AutomorphismGroup   {\AlgebraicClosure  {k}   (x,  y)}
{\AlgebraicClosure   {k}}$  and   construct  $k'$   for  $G=\AutomorphismGroup
{\AlgebraicClosure {k} (x, y)} {\AlgebraicClosure  {k}} / \CyclicGroup {2}$ as
explained above. We show how to apply this method in the following example.

\begin{example}\label{ex:computingAutOverGivenConstantField}
  We  consider  $F:=\GF{7}(x,y)$,  $y^2=  x^5+x^3+x$,  i.e.\  we  examine  the
  curve  from  example~\ref{ex:computingAutOverAlgClosure} over  $\GF{7}$.  We
  already  know,   that  $\AutomorphismGroup  {\AlgebraicClosure{\GF{7}}(x,y)}
  {\AlgebraicClosure{\GF{7}}} \cong \DihedralGroup{6} \times \CyclicGroup{2}$.
  Thus   we   set   $G:=   \DihedralGroup{6}$.  To   find   out,   for   which
  extension   \mbox{$k   \supseteq   \GF{7}$}  we   have   $\AutomorphismGroup
  {k(x,y)}{k}     =    \AutomorphismGroup     {\AlgebraicClosure{\GF{7}}(x,y)}
  {\AlgebraicClosure{\GF{7}}} \cong \DihedralGroup{6} \times \CyclicGroup{2}$,
  we have a closer look at the proof\/\footnote{i.e.\ the computations proving
  that  $\GF{7}(x,y)$ is  indeed of  the specified  type.} that  $\GF{7}(x,y)$
  is  of   type  $\typeOfField{\DihedralGroup{6}}{\AlgebraicClosure{\GF{7}}}$.
  As    seen   in    example~\ref{ex:computingAutOverAlgClosure},   we    have
  $k(x,y)=k(t,u)$,  $u^2=t^6+1$ and  the  automorphism group  is generated  by
  $\varphi:  t \mapsto  t$, $u  \mapsto  -u$, $\psi:  t \mapsto  \zeta t$,  $u
  \mapsto u$  and $\sigma: t \mapsto  \frac 1t$, $u \mapsto  \frac {u} {t^3}$,
  with a primitive $6$-th root of unity $\zeta$. As $3$ is such a $6$-th root,
  we may set $\zeta:=3$. Thus, our  automorphism are defined over the smallest
  extension $k \supseteq \GF{7}$ s.th.\ $t,u \in k(x,y)$.
  
  Hence, to compute $k$ we have to  examine $t$ and $u$ more closely. They can
  be  computed  from  $x$  and  $y$ using  the  coefficients  $\alpha_i$  from
  theorem~\ref{thm:relationOfGeneratorsOfHyperellipticFields}.  Therefore, $k$
  is the  smallest field  s.th.\ $\alpha_i \in  k$. Solving  the corresponding
  equations  and   inequalities,  we  get  that   $x\in  k[t]$,  $\alpha_0=1$,
  $\alpha_1^3=i$ with $i^2=-1$  is a possible solution.  Furthermore, there is
  no solution  over $\GF{7}$.  Thus, $t,u  \in \GF{49}(x,y)  := \GF{7}(x,y,i)$
  which implies  that $k:=\GF{49}$ is  the smallest constant field  s.th. $Fk=k(x,y)$
  has the automorphism group $\DihedralGroup{6}\times \CyclicGroup{2}$.
\end{example}

\section{Computational Aspects}\label{sec:implementation}

The  author   implemented  algorithm~\ref{alg:compAutWrtAlgClosedK}   for  the
computer algebra systems MuPAD  (\cite{MuPAD}) and Singular (\cite{SINGULAR}).
The  Gr\"obner  basis  steps  are implemented  for  Singular,  while  anything
else---i.e.\  Brandt's  normal  forms,  computing $N$,  substitution  and  the
comparing of coefficients---is programmed for MuPAD. Both parts of the program
are combined  using shell scripts.  It was  decided to separate  the Gr\"obner
basis steps from the rest of the  computation, since on the one hand, Singular
has one  of the most efficient  Gr\"obner basis implementations. On  the other
hand,  Singular is  restricted to  characteristic $p\le  32003$, which  is too
small for many fields of cryptographic relevance.

As a proof of  concept, the implementation is not optimized  for speed at all.
Therefore,  a speedup  by a  factor  of at  least  $10$ ought  to be  possible
using  a  ``proper''  implementation.  Nevertheless,  the  examples  given  in
table~\ref{tab:runningTimeExamples}  suggest  that  even  this  implementation
computes  the  automorphism  group  $\AutomorphismGroup  {\AlgebraicClosure{k}
(x,y)} {\AlgebraicClosure{k}}$  of an  arbitrary hyperelliptic  function field
very efficiently.
The      performance     
seems to depend neither   on  the  size
of   the   constant  field,   nor   on   the  order   of   $\AutomorphismGroup
{\AlgebraicClosure{k} (x,y)}  {\AlgebraicClosure{k}}$. Even  though increasing
the genus  increases the size  of the  systems of polynomials---the  number of
both the  polynomials and the  parameters increase  linear with $g$  for types
like  $\typeOfField{\CyclicGroup{2}}{\AlgebraicClosure{k}}$---,  the  examples
indicate that even for fields of  genus $4$ and higher, the automorphism group
computations are quite fast.

\begin{table}\label{tab:runningTimeExamples}
  \begin{tabular}{llr|r}
    $k$ 
       & Defining Equation 
       & $\cardinality{\AutomorphismGroup 
         {\AlgebraicClosure{k} (x,y)} {\AlgebraicClosure{k}}}$ 
       & seconds

    \\\hline 
    {$g=2$:}&&&
    \\
     $\GF{9491}$ &
       $y^2=x^5 - 4608x + 1124$ &
       2 &
       12.6
    \\
     $\GF{10223}$ &
       $y^2=x^6 - 4x^4 - 4x^2 + 1$ &
       4 &
       52.5
    \\
     $\GF{10711}$ &
       $y^2=x^6 + 394x^3 - 3378$ &
       12 &
       23.3
    \\
     $\GF{3}$ &
       $y^2=x^6 + x^4 + x^2 + 1$ &
       24 &
       9.8
    \\
       $\GF{3}$ &
       $y^2=x^6 + x^4 + x^2 + 1$ &
       48 &
       9.2
    \\ 
     $\GF{5}$ &
       $y^2=x^5 + 4x$ &
       240 &
       22.7
    \\\hline
    $g=3$: &&&
    \\
     $\GF{11}$ &
       $y^2=x^7 + 6x^6 + 5x^4 + 4x^3 + x + 3$ &
       2 &
       67.0
    \\
     $\GF{3}$ &
       $y^2=x^8 + x^7 + 2x^5 + 2 x + 2$ &
       8 &
       30.3
    \\
     $\GF{7}$ &
       $y^2=x^7 + 6 x^4 + 4x^3 + x^2 + 2$ &
       42 &
       67.8
    \\\hline
    $g=4$: &&&
    \\
     $\GF{5}$ &
       $y^2=x^{10} + x^8 + 3x^6 + 4x^2 + 4$ &
       4 &
       81.8
    \\ 
     $\GF{3}$ &
       $y^2=x^9 + 2x^7 + 2x^3 + 2x$ &
       8 & % Das semidirekte Produkt kann's nicht sein, da keine C_3 gefunden wurde!
       46.6
    \\ 
%     $\GF{3}$ &
%       $y^2=x^9 +2x^3 + x +2$ &
%       36 &
%       27.7
%    \\ 
%     $\GF{3}$ &
%       $y^2=x^9 + 2x$ &
%       1440 &
%       34.1
%    \\ 
  \end{tabular}\par\mbox{}
  \caption{Time to compute ${\AutomorphismGroup 
         {\AlgebraicClosure{k} (x,y)} {\AlgebraicClosure{k}}}$
	 on an Intel$^{\text{\textregistered}}$ Ce\-le\-ron$^{\text{\textregistered}}$, 
	 \unit[1.7]{GHz}, ordered by genus and $\cardinality{\AutomorphismGroup 
         {\AlgebraicClosure{k} (x,y)} {\AlgebraicClosure{k}}}$ 
	 }
\end{table}

Let  us  discuss  the  cryptographic application,  briefly.  As  explained  in
the  introduction, the initial goal  was to  provide  an algorithm  to
check,  whether  a  given  hyperelliptic  curve promises  to  yield  a  secure
Jacobian,  i.e.\ whether  it  is  worthwhile to  apply  more expensive  algorithms
to  check  a given  curve  for  security.  Because  of the  attacks  mentioned
in   the  introduction,   secure   curves  have   small  automorphism   groups
$\AutomorphismGroup {k(x,y)} {k}$. Since $\AutomorphismGroup {k(x,y)} {k}
\le  \AutomorphismGroup  {\AlgebraicClosure{k}(x,y)}  {\AlgebraicClosure{k}}$,
algorithm~\ref{alg:compAutWrtAlgClosedK}  can be used to assure 
this property.
The  timings  of table~\ref{tab:runningTimeExamples}  also  apply  to the  set
of  relevant  curves,  as secure  curves  are  of  genus  $\le 4$  because  of
the Adleman-DeMarrais-Huang attack  (\cite{adlemanDemarraisHuang:1994}) and as
characteristic of the constant field and the size of the automorphism group do
not seem to influence the running time. 

Even though a small automorphism group
is necessary for a secure curve, it  is not obvious, how much information concerning
security can be deduced from knowing the automorphism  group. A discussion
of this topic can be found in~\cite{goeb:2003}.
 
The   methods  described   in  section~\ref{subsec:compAut:arbitraryConstants}
were     not     implemented.     Nevertheless,     we     will     try     to
compare    algorithm~\ref{alg:compAutWrtAlgClosedK}    to   Michael    Stoll's
\texttt{Au\-to\-mor\-phism\-Group}    function     (cf.\    \cite{stoll:2001})
in    some    examples.     To    do    so,  we choose  the
smallest    field   $k$    of    the   given    characteristic,   for    which
$\AutomorphismGroup  {k(x,y)} {k}  = \AutomorphismGroup  {\AlgebraicClosure{k}
(x,y)}   {\AlgebraicClosure{k}}$ holds,  in each example.   Then,   $\AutomorphismGroup
{\AlgebraicClosure{k}   (x,y)}  {\AlgebraicClosure{k}}$   is  computed   using
algorithm~\ref{alg:compAutWrtAlgClosedK},  while  Stoll's  method is  used  to
compute $\AutomorphismGroup {k(x,y)} {k}$. The running times for some examples
are given in table~\ref{tab:runningTimeComparisonToStoll}.

\begin{table}\label{tab:runningTimeComparisonToStoll}
  \begin{tabular}{lll|ll}
    \multicolumn{3}{l|}{Function Field}&\multicolumn{2}{l}{Running Time} \\
    $k$ & Defining Equation & $\cardinality {\AutomorphismGroup {k(x,y)} {k}}$ & Stoll & G\"ob \\
    \hline
     $\GF{3^6}$ &
       $y^2=x^9 + 2x^3 + x +2$ &
       36 &
       2.9 &
       27.7
    \\ 
     $\GF{5^2}$ &
       $y^2=x^5 + 4x$ &
       240 &
       18.1 &
       22.7
    \\
      $\GF{7^2}$ &
      $y^2=x^7+6x$ & 
      672 &
      228.3 &
      61.0 % PGL_2(7)
    \\ 
     $\GF{3^4}$ &
       $y^2=x^9 + 2x$ &
       1440 &
       1347.3 &
       34.1
    \\ 
     $\GF{11^2}$ &
       $y^2=x^{11} +10x$ &
       2640 &
       5625.1 &
       90.3
       5
  \end{tabular}\par\mbox{}
  \caption{Running time comparison between Michael Stoll's algorithm
           and algorithm~\ref{alg:compAutWrtAlgClosedK}, timings in
           seconds on an Intel$^{\text{\textregistered}}$~Celeron$^{\text{\textregistered}}$, 
	   \unit[1.7]{GHz}}
\end{table}

From  these examples,  Stoll's  algorithm seems  to be  quite  fast for  small
automorphism groups,  while it is very  slow for large ones.  As stated above,
our implementation does  not seem to be  influenced by the group  size at all.
Thus,  if  you are  quite  sure  that the  field  you  are investigating  only
has  a small  automorphism group,  Stoll's  algorithm ought  to be  preferred.
Even  though  the  majority  of  hyperelliptic function  fields  has  a  small
automorphism group, the remaining fields do  not seem to be suited for Stoll's
algorithm. Hence, in  order to compute the automorphism group  of an arbitrary
hyperelliptic function field, it might be  sensible to use the algorithms from
section~\ref{subsec:compAut:arbitraryConstants} as  those at least seem  to be
more predictable  w.r.t.\ performance. Furthermore, Stoll's  algorithm returns
every  single  automorphism,  while  the  methods  presented  here,  give  the
structure as well  as the generators of the automorphism  group. Thus, it also
depends on the application, which of the algorithms ought to be used.

%\bibliographystyle{alpha}
%\bibliography{hyperell,fundamentals,groups} 
\newcommand{\etalchar}[1]{$^{#1}$}

\end{document}